\input amstex 
\documentstyle{amsppt}
\magnification 1200
\topmatter
\title
Exactness of martingale approximation and the central limit theorem 
\endtitle
\author
Dalibor Voln\'y
\endauthor
\affil
Universit\'e de Rouen 
\endaffil
\address
D\'epartement de Math\'ematiques, 
Universit\'e de Rouen, 
76801 Saint Etienne du Rouvray,
France
\endaddress
\email{Dalibor.Volny\@univ-rouen.fr}
\endemail
\abstract
Let $(X_i)$ be a Markov chain with kernel $Q$, $f$ an $L^2$ function on its state space. 
If $Q$ is a normal operator and $f = (I-Q)^{1/2}g$ (which is equivalent to the convergence of
$\sum_{n=1}^\infty \frac{\sum_{k=0}^{n-1}Q^kf}{n^{3/2}}$ in $L^2$), by Derriennic and Lin \cite{D-L}
we have the central limit theorem. By \cite{M-Wu} and \cite{Wu-Wo} the CLT is implied by the convergence of 
$\sum_{n=1}^\infty \frac{\|\sum_{k=0}^{n-1}Q^kf\|_2}{n^{3/2}}$ and by $\|\sum_{k=0}^{n-1}Q^kf\|_2 =
o(\sqrt n/\log^q n)$, $q>1$. We shall show that if $Q$ is not normal or if the conditions of Maxwell and Woodroofe
or Wu and Woodroofe are weakened by $\sum_{n=1}^\infty c_n\frac{\|\sum_{k=0}^{n-1}Q^kf\|_2}{n^{3/2}}<\infty$ for some sequence
$c_n\searrow 0$ or by $\|\sum_{k=0}^{n-1}Q^kf\|_2 = O(\sqrt n/\log n)$, the CLT need not hold.
\endabstract
\endtopmatter
\document
\subheading{1. Introduction} Let $(\Omega, \Cal A,\mu)$ be a probability space with a bijective, bimeasurable and
measure preserving transformation $T$. For a measurable function $f$ on $\omega$, $(f\circ T^i)_i$ is a (strictly)
stationary process and reciprocally, any (strictly) stationary process can be represented in this way.

Billingsley and Ibragimov (cf\. \cite{B}, 
\cite{I}) have proved that if $(m\circ T^i)$ is a martingale difference sequence with $m\in L^2$ and $\mu$ is ergodic 
(i.e\. all sets $A\in\Cal A$ for which $A=T^{-1}A$ it is $\mu(A)=0$ or $\mu(A)=0$) then $\frac1{\|m\|_2\sqrt n} S_n(m)$ 
converge in law to the standard normal law $\Cal N(0,1)$.
Since the publication of Gordin's contribution \cite{G}, a special attention has been given to proving limit theorems
via approximations by martingales. An important part of such results concern Markov chains, cf\. e.g\. \cite{G-L}, \cite{K-V}, \cite{W}, \cite{D-L}, \cite{Wu-W}.

Let $(S,\Cal B,\mu)$ be a probability space, $(\xi_i)$ a homogeneous and ergodic Markov chain with state space $S$,
transition operator $Q$, and stationary distribution $\mu$. For a measurable function $f$ on $S$, $(f(\xi_i))$ is then a
stationary random process; we shall study the central limit theorem for 
$$
  S_n(f) = \sum_{i=0}^{n-1} f(\xi_i)
$$
where $f\in L_0^2(\mu)$, i.e\. is square integrable and has zero mean. Gordin and Lif\v sic (\cite{G-L}) showed that if $f$
is a solution of the equation
$$
  f = g - Qg
$$
with $g\in L^2$ then a martingale approximation giving the CLT exists. The result was extended to reversible operators $Q$
and $f$ satisfying
$$
  f = (I-Q)^{1/2}g \tag1
$$
with $g\in L^2$ by Kipnis and Varadhan in \cite{K-V}, then for normal operators $Q$ and $f$ satisfying (1) by Derriennick and Lin
in \cite{D-L}. As noticed by Gordin and Holzmann (\cite{G-H}), (1) is equivalent to the convergence of
$$
  \sum_{n=1}^\infty \frac{\sum_{k=0}^{n-1}Q^kf}{n^{3/2}}\quad\text{in}\quad L^2. \tag2
$$
Maxwell and Woodroofe have shown in \cite{M-W} that if
$$
  \sum_{n=1}^\infty \frac{\sum_{k=0}^{n-1}\|Q^kf\|_2}{n^{3/2}} < \infty \tag3
$$
(without any other assumptions on the Markov operator $Q$) then the CLT takes place. 

Remark that any stationary process can be represented by a Markov chain
(cf\. \cite{Wu-W}). The central limit theorem of Maxwell and Woodroofe can thus be expressed in the following way:

Let $(\Omega, \Cal A, P, T)$ be a probability space with a bimeasurable and measure preserving bijective transformation $T : \Omega
\to \Omega$, $\Cal F_i$ an increasing filtration with $T^{-1}\Cal F_i = \Cal F_{i+1}$, $f$ is a square integrable and zero mean function
on $\Omega$, $\Cal F_0$-measurable. We denote
$$
  S_n(f) = \sum_{i=0}^{n-1} f\circ T^i. 
$$
(2) then becomes
$$
  \sum_{n=1}^\infty \frac{E(S_n(f)\,|\,\Cal F_0)}{n^{3/2}} \tag2'
$$
and (3) becomes
$$
  \sum_{n=1}^\infty \frac{\|E(S_n(f)\,|\,\Cal F_0)\|_2}{n^{3/2}} < \infty. \tag3'
$$ 
In \cite{Wu-W}, Wu and Woodroofe have shown that if
$$
  \|E(S_n(f)\,|\,\Cal F_0)\|_2 = o\Big(\frac{\sqrt n}{n^q}\Big) \tag4
$$
for some $q>1$ then the CLT takes place. 

In \cite{Vo} and \cite{Kl-Vo1}, nonadapted versions of Maxwell-Woodroofe aproximations (3') and have been found. 

In the present paper we will deal with exactness of the central limit theorems of Derriennic 
and Lin, Maxwell and Woodroofe, and of Wu and Woodroofe. 
First, we show that the central limit theorem of Derriennic and Lin cannot be extended to non normal operators $Q$. 

\proclaim{Theorem 1} There exists a process $(f\circ T^i)$ such that the series
$$
  \sum_{n=1}^\infty \frac{E(S_n(f)\,|\,\Cal F_0)}{n^{3/2}} \tag2'
$$
converges in $L^2$,
but for two different subsequences $(n_k')$, $(n_k'')$, the distributions of $S_{n_k'}/\sigma_{n_k'}$ and
$S_{n_k''}/\sigma_{n_k''}$ converge to different limits.
\endproclaim

Then we show that in the central limit theorems of Maxwell and Woodroofe and of Wu and Woodroofe, the rate of convergence
of $\|E(S_n(f)|\Cal F_0)\|_2$ towards 0 is practically optimal.

\proclaim{Theorem 2} For any sequence of positive reals $c_n\to 0$ there exists a process $(f\circ T^i)$ such that 
$$
  \sum_{n=1}^\infty c_n\frac{\|E(S_n(f)|\Cal F_0)\|_2}{n^{3/2}} < \infty \tag5
$$
but for two different subsequences $(n_k')$, $(n_k'')$, the distributions of $S_{n_k'}/\sigma_{n_k'}$ and
$S_{n_k''}/\sigma_{n_k''}$ converge to different limits.
\endproclaim

In \cite{Pe-U}, Peligrad and Utev have shown that under the same assumptions there exists an $f$
such the sequence of $S_n(f)/\sqrt n$ is not stochastically bounded. (Remark that in the same 
paper the authors have proved that (3') implies also the weak invariance principle.)
 
\proclaim{Theorem 3} There exists a process $(f\circ T^i)$ such that 
$$
  \|E(S_n(f)|\Cal F_0)\|_2 = O\Big(\frac{\sqrt n}{\log n}\Big),  \tag6 
$$
but for two different subsequences $(n_k')$, $(n_k'')$, the distributions of $S_{n_k'}/\sigma_{n_k'}$ and
$S_{n_k''}/\sigma_{n_k''}$ converge to different limits.
\endproclaim

From the construction it follows that in Theorems 1-3, the variances $\sigma_n^2$ of $S_n(f)$ grow faster than linearly.
It thus remains an open problem whether with a supplementary assumption $\sigma_n^2/n\to const.$ the CLT would hold.
As shown by a couter example in \cite{Kl-Vo2}, this assumption is not sufficient for
$q\leq 1/2$, the only exponents to consider are thus $1/2<q<1$. 

It also remains an open question whether the CLT would hold for $f\in L^{2+\delta}$ for some $\delta>0$.

\subheading{2. Proof}

We give one proof which will treat all three theorems.

In all of the text, $\log$ will denote the dyadic logarithm.

For $k=1,2,\dots$ let $n_k = 2^k$, $e_k$ be random variables with
$$
  \|e_k\|_2 = a_k/k,\quad 0\leq a_k\leq 1,\quad \sum_{k=1}^\infty a_k/k =\infty,
$$ 
such that for each $k$, 
$U^ie_k$ are independent, and if $i\neq j$ then $U^ie_{k'}$ and $U^je_{k''}$ are orthogonal. For $k'\neq k''$ the 
$e_{k'}$, $e_{k''}$ are not orthogonal but we suppose that for all $1\leq k',k''$ it is $E(e_{k'}e_{k''})\geq 0$ and
$$
  \Big\| \sum_{k=1}^n e_k \Big\|_2 \nearrow \infty \quad\text{as}\quad n\to \infty.
$$
Let
$$
  f = \sum_{k=1}^\infty \frac1{n_k} \sum_{i=0}^{n_k-1} U^{-i}e_k.
$$
We have $\|f\|_2 \leq \sum_{k=1}^\infty \|e_k\|_2/\sqrt{n_k} <\infty$ due to the exponential growth of the $n_k$s.

For a positive integer $N$ we have
$$
  S_N(f) = \sum_{k=1}^\infty \sum_{j=0}^{N-1} \sum_{i=0}^{n_k-1} \frac1{n_k} U^{j-i}e_k = 
  S'_N(f) + S''_N(f)
$$
where
$$
  S'_N(f) = S_N(f) - E(S_N(f)|\Cal F_0) = \sum_{k=1}^\infty \sum_{j=0}^{N-1} 
  \sum_{i=0}^{(j\wedge n_k)-1} 
  \frac1{n_k} U^{j-i}e_k
$$
($j\wedge n_k = \min\{j, n_k\}$) and
$$
  S''_N(f) = E(S_N(f)|\Cal F_0) = \sum_{k=1}^\infty \sum_{j=0}^{N-1} \sum_{i=j}^{n_k-1} 
  \frac1{n_k} U^{j-i}e_k.
$$

We will study the asymptotic behaviour of $S''_N(f) = E(S_N(f)|\Cal F_0)$ and $S'_N(f) = S_N(f) - E(S_N(f)|\Cal F_0)$
separately. In the first case
we will show that the series
$$
  \sum_{n=1}^\infty \frac{E(S_n(f)\,|\,\Cal F_0)}{n^{3/2}} 
$$
converges in $L^2$, for a suitable choice of $a_k$ we shall have 
$$
  \|E(S_n(f)|\Cal F_0)\|_2 = O\Big(\frac{\sqrt n}{\log n}\Big) 
$$
and for any sequence of positive reals $c_n\to 0$ the $a_k$ can be chosen so that 
$$
  \sum_{n=1}^\infty c_n\frac{\|E(S_n(f)|\Cal F_0)\|_2}{n^{3/2}} < \infty. 
$$

In the second case we will show that the assumption $\sum_{k=1}^\infty a_k/k =\infty$
allows us to define the $e_k$ so that for two different subsequences $(n_k')$, $(n_k'')$, the distributions 
of $S_{n_k'}/\sigma_{n_k'}$ and $S_{n_k''}/\sigma_{n_k''}$ converge to different limits.

Eventually we prove existence of a dynamical system on which the process can be defined.
\bigskip
\centerline{\it 1. Asymptotics of $S''_N(f) = E(S_N(f)|\Cal F_0)$.}
\medskip
For $N\leq n_k$ we have
$$
  \sum_{j=0}^{N-1} \sum_{i=j}^{n_k-1} U^{j-i}e_k = \sum_{i=0}^{n_k-N} N U^{-i}e_k + 
  \sum_{i=n_k-N+1}^{n_k-1}
  (n_k-i) U^{-i}e_k
$$
and for $N>n_k$ we have
$$
  \sum_{j=0}^{N-1} \sum_{i=j}^{n_k-1} U^{j-i}e_k = \sum_{j=0}^{n_k-1} (n_k-j) U^{-j}e_k,
$$
hence
$$\multline
  S''_N(f) = \sum_{k\geq 1: n_k< N} \sum_{j=0}^{n_k-1} \frac{n_k-j}{n_k} U^{-j}e_k + \\
  + \sum_{k\geq 1: n_k\geq N} \Big[\sum_{j=0}^{n_k-N} \frac{N}{n_k} U^{-j}e_k +
  \sum_{j=n_k-N+1}^{n_k-1} \frac{n_k-j}{n_k} U^{-j}e_k \Big].
  \endmultline \tag7
$$
We will prove that (2') is satisfied. For this, it is sufficient to show that
$$
  \sum_{N=1}^\infty \sum_{k\geq 1: n_k< N} \sum_{j=0}^{n_k-1} 
  \frac1{N^{3/2}}\frac{n_k-j}{n_k} U^{-j}e_k
$$
and
$$
  \sum_{N=1}^\infty \sum_{k\geq 1: n_k\geq N} \frac1{N^{3/2}} \Big[\sum_{j=0}^{n_k-N} 
  \frac{N}{n_k} U^{-j}e_k +
  \sum_{j=n_k-N+1}^{n_k-1} \frac{n_k-j}{n_k} U^{-j}e_k \Big]
$$
converge in $L^2$.

Recall that $n_k=2^k$. For the first sum we have
$$\gather
  \Big\|\sum_{N=1}^\infty \sum_{k\geq 1: 2^k< N} \sum_{j=0}^{2^k-1} 
  \frac1{N^{3/2}}\frac{2^k-j}{2^k} U^{-j}e_k \Big\|_2^2 = \\
  \Big\|\sum_{k=1}^\infty \Big(\sum_{N=2^k+1}^\infty \frac1{N^{3/2}}\Big)
  \sum_{j=0}^{2^k-1} \frac{2^k-j}{2^k} U^{-j}e_k \Big\|_2^2 \leq \\
  \Big\| c\sum_{k=1}^\infty \sum_{j=0}^{2^k-1} \frac{2^k-j}{2^{3k/2}} U^{-j}e_k \Big\|_2^2 
  \leq \\
  \Big\| c\sum_{j=0}^\infty \sum_{k>\log (j+1)} \frac{2^k-j}{2^{3k/2}} U^{-j}e_k \Big\|_2^2 
  \leq \\
  \sum_{j=0}^\infty \Big\| c\sum_{k>\log (j+1)} \frac1{2^{k/2}} U^{-j}e_k \Big\|_2^2 \leq 
  C \sum_{j=0}^\infty \Big(\frac1{\sqrt{(j+1)}[1\vee\log (j+1)]}\Big)^2 < \infty
  \endgather
$$
where $0<c,C<\infty$.
For the second sum we have
$$\gather
  \Big\|\sum_{N=1}^\infty \sum_{k\geq 1: 2^k\geq N} \frac1{N^{3/2}} \sum_{j=0}^{2^k-N} 
  \frac{N}{2^k} U^{-j}e_k \Big\|_2^2 \leq \\
  \Big\|\sum_{N=1}^\infty \sum_{j=0}^\infty \sum_{k\geq \log (N+j)} \frac{1}{N^{1/2}2^k} 
  U^{-j}e_k \Big\|_2^2 \leq
  \Big\| \sum_{j=0}^\infty \sum_{k\geq \log(j+1)} \Big(\sum_{N=1}^{2^k} 
  \frac{1}{N^{1/2}}\Big) \frac{1}{2^k} U^{-j}e_k \Big\|_2^2 \leq \\
  \sum_{j=0}^\infty \Big\| c \sum_{k\geq \log(j+1)} \frac1{2^{k/2}} U^{-j}e_k \Big\|_2^2 
  \leq C
  \sum_{j=0}^\infty \Big(\frac1{\sqrt{(j+1)}[1\vee\log (j+1)]}\Big)^2 < \infty
\endgather
$$ 
and
$$\gather
  \Big\| \sum_{N=1}^\infty \sum_{k\geq 1: 2^k\geq N} \frac1{N^{3/2}} 
  \sum_{j=2^k-N+1}^{2^k-1} \frac{2^k-j}{2^k} U^{-j}e_k \Big\|_2^2 \leq \\
  \sum_{j=0}^\infty \Big\|\sum_{k\geq \log(j+1)} 
  \sum_{N=2^k-j+1}^{2^k}\frac1{N^{3/2}}\frac{2^k-j}{2^k} U^{-j}e_k \Big\|_2^2 \leq \\
  \sum_{j=0}^\infty \Big\|\sum_{k\geq \log(j+1)} 
  \Big(\sum_{N=2^k-j+1}^{2^k}\frac{1}{N^{1/2}}\Big) \frac{1}{2^k} U^{-j}e_k \Big\|_2^2 \leq \\
  C \sum_{j=0}^\infty \Big(\frac1{\sqrt{(j+1)}[\log (j+1)]}\Big)^2 < \infty
\endgather
$$
where $0<c,C<\infty$. This finishes the proof of (2').
\medskip

We have
$$\gathered
  (1/\sqrt 6)\|e_k\|_2\sqrt{n_k}\leq \Big\|\sum_{j=0}^{n_k-1} \frac{n_k-j}{n_k} 
  U^{-j}e_k\Big\|_2 \leq \|e_k\|_2\sqrt{n_k},
  \quad n_k<N,\\
  \Big\|\sum_{j=0}^{n_k-N} \frac{N}{n_k} U^{-j}e_k + \sum_{j=n_k-N+1}^{n_k-1} 
  \frac{n_k-j}{n_k} U^{-j}e_k \Big\|_2
  \leq \frac{N}{\sqrt{n_k}} \|e_k\|_2, \quad n_k\geq N.
  \endgathered \tag8
$$

Recall that by $[x]$ we denote the integer part of $x$.
Because $n_k = 2^k$ grow exponentially fast, there exists a constant $0<c<\infty$ not 
depending on $N$ such that 
$$
  \sum_{k\geq 1: n_k\geq N} (N/\sqrt{n_k})\|e_k\|_2 \leq 
  c\|e_{[\log N]}\|_2,\quad \sum_{k\geq 1: n_k< N} \|e_k\|_2\sqrt{n_k} \leq c \sqrt N
  \|e_{[\log N]}\|_2.
$$ 
Using (7) and (8) we deduce that for some constants $c',c''>0$ we have
$$
  c'N \|e_{[\log N]}\|_2^2< \|E(S_N(f)|\Cal F_0)\|_2^2 <c''N \|e_{[\log N]}\|_2^2.
$$
Because $\|e_k\|_2 = a_k/k$,
$$
  c'\frac{Na^2_{[\log N]}}{\log^2 N} < \|E(S_N(f)|\Cal F_0)\|_2^2 <c''
  \frac{Na^2_{[\log N]}}{\log^2 N}.\tag9
$$
For $a_k\equiv 1$ we thus have 
$$
  \|E(S_n(f)|\Cal F_0)\|_2 = O\Big(\frac{\sqrt n}{\log n}\Big).
$$

From (9) we deduce that the series $\sum_{n=1}^\infty n^{-3/2}\|E(S_n(f)|\Cal F_0)\|_2$ converges if and only if
$\sum_{n=1}^\infty a_{[\log n]}/(n[\log n])$ converges; because
$$
  \frac{a_k}{2k} \leq \sum_{j=0}^{2^k-1} \frac{a_k}{k(2^k+j)} \leq \frac{a_k}{k}
$$
this is equivalent to the convergence of
$\sum_{n=1}^\infty a_n/n$.

Let $c_n$ be positive real numbers, $c_n\to 0$; we can choose the $a_n$ so that
$\sum_{n=1}^\infty a_nc_n/n<\infty$, that means
$$
  \sum_{n=1}^\infty c_n\frac{\|E(S_n(f)|\Cal F_0)\|_2}{n^{3/2}} \approx 
  \sum_{n=1}^\infty \frac{c_na_{[\log n]}}{n[\log n]} <\infty
$$
but $\sum_{n=1}^\infty a_n/n = \infty$, i.e\. .
$$
  \sum_{n=1}^\infty \frac{\|E(S_n(f)|\Cal F_0)\|_2}{n^{3/2}} = \infty.
$$


\bigskip
\centerline{\it 2. Asymptotics of $S'_N(f) = S_N(f) - E(S_N(f)|\Cal F_0)$.}
\medskip

Notice that in the preceding section, no hypothesis on dependence of the $e_k$ was needed. Now, we shall suppose that the
sequence of $a_k$ is fixed and we choose the $e_k$ so that for two different subsequences $(n_k')$, $(n_k'')$, the distributions 
of $S_{n_k'}/\sigma_{n_k'}$ and $S_{n_k''}/\sigma_{n_k''}$ converge to different limits.

For $N\leq n_k$ we have
$$
  \sum_{j=0}^{N-1} \sum_{i=0}^{(j\wedge n_k)-1} U^{j-i}e_k = \sum_{j=1}^{N-1} (N-j)U^je_k 
  \tag{10}
$$
and for $N>n_k$ we have
$$
  \sum_{j=0}^{N-1} \sum_{i=0}^{(j\wedge n_k)-1} U^{j-i}e_k = 
  \sum_{j=1}^{N-n_k} n_kU^je_k + \sum_{j=N-n_k+1}^{N-1} (N-j)U^je_k. \tag{11}
$$
For all $k\geq 1$ we have $PlU^je_k = 0$ if $j\neq l$, $PlU^le_k = U^le_k$.
For $l\geq N$ and $l\leq 0$ we thus have $P_lS_N(f) = 0$ and for $1\leq l\leq N-1$ we, using (10) and (11), deduce
$$
  P_lS_N(f) = \sum_{k\geq 1: n_k\leq N-l} U^le_k + \sum_{k\geq 1: n_k\geq N+1-l} \frac{N-l}{n_k} U^le_k. \tag{12}
$$
\comment
For any $1>\epsilon>0$ (arbitrarilly small) and $1\leq l\leq N-\epsilon N=N(1-\epsilon)$,  
$$
  \sum_{k\geq 1: n_k\leq N-l} U^le_k = \sum_{k\geq 1: n_k\leq \epsilon N} U^le_k + 
  \sum_{k\geq 1: \epsilon N<n_k\leq N} U^le_k.
$$
\endcomment
Recall that $[x]$ denotes the integer part of $x$. We have
$$
  S'_N(f) = \sum_{l=1}^{N-1} P_lS_N(f) = \sum_{l=1}^{[N(1-\epsilon)]} P_lS_N(f) +
  \sum_{l=[N(1-\epsilon)]+1}^{N} P_lS_N(f)
$$
where
$$\multline
  \sum_{l=1}^{[N(1-\epsilon)]} P_lS_N(f) = \\
  \sum_{l=1}^{[N(1-\epsilon)]} U^l \sum_{k\geq 1: n_k\leq N-l} e_k +
  \sum_{l=1}^{[N(1-\epsilon)]} U^l \sum_{k\geq 1: n_k\geq N+1-l} \frac{N-l}{n_k} e_k = \\
  = \sum_{l=1}^{[N(1-\epsilon)]} U^l \sum_{k\geq 1: n_k\leq \epsilon N} e_k +
  \sum_{l=1}^{[N(1-\epsilon)]} U^l \sum_{k\geq 1: \epsilon N<n_k\leq N-l} e_k + \\
  + \sum_{l=1}^{[N(1-\epsilon)]} U^l \sum_{k\geq 1: n_k\geq N+1-l} \frac{N-l}{n_k} e_k.
  \endmultline
$$
Because $n_k=2^k$, 
$$
  \Big\|\sum_{k\geq 1: n_k\geq N+1-l} \frac{N-l}{n_k} U^le_k\Big\|_2 \leq 2\|e_{[\log (N-l)]}\|_2 \leq 2/\log (N-l); \tag{13}
$$ 
$\epsilon N <n_k\leq N$ if and only if $\log N + \log\epsilon <k\leq \log N$. We thus deduce that for $\epsilon>0$ fixed and
$b(N) = \|\sum_{k=1}^{\log N} e_k\|_2\nearrow \infty$ ,
$$
  \lim_{N\to\infty} \frac1{b(N)\sqrt N} \Big\|\sum_{\epsilon N<n_k\leq N} e_k\Big\|_2 = 0
$$
and
$$
  \lim_{N\to\infty} \frac1{b(N)\sqrt N} \Big\|\sum_{l=1}^{[N(1-\epsilon)]} P_lS_N(f) - \sum_{l=1}^{[N(1-\epsilon)]} U^l
  \sum_{k\geq 1: n_k\leq \epsilon N} e_k\Big\|_2 = 0.
$$
For all $[N(1-\epsilon)]+1\leq l\leq N-1$ we have, by (12) and (13),
$\|P_lS_N(f)\|_2 \leq b(N-l) + 2/\log (N-l)$ hence 
$$
  \lim_{\epsilon\searrow 0}\lim_{N\to\infty} \frac1{b(N)\sqrt N} \Big\|\sum_{l=[N(1-\epsilon)]+1}^{N-1} P_lS_N(f) - 
  \sum_{l=[N(1-\epsilon)]+1}^{N-1} \sum_{k\geq 1: n_k\leq \epsilon N} U^le_k\Big\|_2 = 0,
$$
therefore
$$
  \lim_{N\to\infty} \frac1{b(N)\sqrt N} \Big\|S'_N(f) - \sum_{l=0}^{N-1} U^l \sum_{k\geq 1: n_k\leq N} e_k\Big\|_2 = 0. \tag{14}
$$

Recall that 
$$
  \|e_k\|_2 = a_k/k,\quad 0\leq a_k\leq 1,\quad \sum_{k=1}^\infty a_k/k =\infty.
$$
Let $N_l$, $l=1,2,\dots$, be an increasing sequence of positive integers such that 
$$
  2^{2^l}-1 < \sum_{k\geq 1: N_{l-1}< n_k\leq N_l} \frac{a_k}{k} < 2^{2^l}+1; \tag{15}
$$
we suppose that for $N_{l-1}< n_k\leq N_l$ the random variables $e_k$ are multiples one of
another and are independent of any $e_j$ with $n_j\leq N_{l-1}$ or $n_j > N_l$. For $l$ odd
we choose $e_k$, $N_{l-1}< n_k\leq N_l$, so that 
$$
  \frac1{b(N_l)\sqrt N_l} \sum_{j=0}^{N_l-1} U^j \Big(\sum_{k\geq 1: N_{l-1}<n_k\leq N_l} e_k\Big)
$$
weakly converge to a symmetrised Poisson distribution and for $l$ even to the standard normal
distribution. We can do so by defining, for $l$ odd and $N_{l-1}< n_k\leq N_l$, $e_k=\pm r_k$ with probabilities
$1/(2N_l)$ and $e_k=0$ with probability $1-1/N_l$, where
$$
  \sum_{k\geq 1: N_{l-1}< n_k\leq N_l} r_k = b(N_l)\sqrt N_l,
$$
for $l$ even we define $e_k$ normally distributed with zero means and variances $r_k^2$, 
$$
  \sum_{k\geq 1: N_{l-1}< n_k\leq N_l} r_k = b(N_l)\sqrt N_l.
$$
By (14) and (15) we then get the convergence to the same laws of $1/(b(N_{2l})\sqrt{N_{2l}})S_{N_{2l}}(f)$
and $1/(b(N_{2l-1})\sqrt{N_{2l-1}})S_{N_{2l-1}}(f)$.
\bigskip
\centerline{\it 3. Existence of $(\Omega,\Cal A,\mu,T)$.}
\medskip

We define, for $l=1,2,\dots$, $A_l = \{-1,0,1\}$ for $l$ odd and $A_l=\Bbb R$ for $l$ even, equipped with the probability measures
$\nu_l$ such that $\nu_l(\{-1\}) = 1/(2N_l) = \nu_l(\{1\})$ and $\nu_l(\{0\}) = 1-1/N_l$ for $l$ odd, $\nu_l = \Cal N(0,1)$ for $l$ even.
For each $l$ we define $\Omega_l = A_l^\Bbb Z$ equipped with the product measure $\mu_l=\nu_l^\Bbb Z$ and with the 
transformation $T_l$ of the left shift. Then we put $\Omega = \underset l=1 \to{\overset \infty \to{\times}} \Omega_l$,
equipe it with the product measure $\mu = \underset l=1 \to{\overset \infty \to{\otimes}} \mu_l$ and the product
transformation $T$. The random variables $e_k$, $N_{l-1}<n_k\leq N_l$, will then be multiples of projections of $\Omega$
onto $A_l$.
\smallskip

This finishes the proof.

\enddocument
\end

\ref \key  \by  \paper \jour \vol  \pages \yr \endref

\enddocument
\end

Let $(\Omega, \Cal A,\mu)$ be a probability space with a bijective, bimeasurable and
measure preserving transformation $T$, $(\Cal F_i)_{i\in\Bbb Z}$ is a filtration, $T^{-1}\Cal F_i = \Cal F_{i+1}$. 
By $f$ we denote an $\Cal F_0$-measurable function with $E(f) =0$ and $E(f^2)<\infty$. We shall suppose that $f$ 
is regular, i.e\. $E(f|\Cal F_{-\infty})=0$ and $E(f|\Cal F_{\infty})=f$. We shall, moreover, suppose that the
process $(\circ T^i)_if$ is adapted, i.e\. $f$ is $\Cal F_0$-measurable.

In \cite{M-Wo} Maxwell and Woodroofe have proved that if the series
$$
  \sum_{n=1}^\infty \frac{\|E(S_n(f)|\Cal F_0)\|_2}{n^{3/2}} \tag1
$$
converges then the (conditional) central limit theorem holds. In \cite{Wu-Wo} Wu and Woodroofe have proved that
for $\sigma_n = \|S_n(f)\|_2$, $\|E(S_n(f)|\Cal F_0)\|_2 =o(\sigma_n)$ is equivalent to an existence of an
array of row-wise stationary martingale difference sequences $(D_{n,k})$, $M_n= \sum_{j=1}^n D_{n,j}$, such that
$\|S_n(f)- M_n\|_2 = o(\sigma_n)$. Because the central limit theorem for martingale difference sequences has been 
largely studied, this approximation gives a nice way for a study of the limit behaviour of the partial sums
$S_n(f)$. In \cite{Wu-Wo}, using this approximation, central limit theorems for processes with nonlinear growth of
variance have been proved.

In \cite{Wu-Wo} Wu and Woodroofe have proved that if, for a $q>1$, 
$$
  \|E(S_n(f)|\Cal F_0)\|_2 = o\Big(\frac{\sqrt n}{\log^qn}\Big) \tag2
$$
then there exists a stationary sequence of martingale differences $(D_j)$ such that $\|S_n(f)- \sum_{j=1}^n D_{j}\|_2 = o(\sigma_n)$ and the central limit theorem thus holds (cf\. \cite{Ha-He, Chapter 5}).

In \cite{Pe-U}, Peligrad and Utev showed that for any sequence of positive reals $c_n\to 0$ there exists an $f$
such that
$$
  \sum_{n=1}^\infty c_n\frac{\|E(S_n(f)|\Cal F_0)\|_2}{n^{3/2}} < \infty \tag3
$$
but the sequence of $S_n(f)/\sqrt n$ is not stochastically bounded. (In the same paper the authors have proved
that (1) implies also the weak invariance principle.)

Let us remark that in \cite{Vo}, a nonadapted version of the Wu-Woodroofe aproximation as well of (1) has been found. 
Using the same method, a nonadapted version of (2) can be given (\cite{K}).
\medskip

In this paper we show that the Wu-Woodroofe aproximation $\|S_n(f)- M_n\|_2 = o(\sigma_n)$ need not imply the CLT
even if the variances $\sigma_n^2$ grow linearly. In particular, we shall prove:

\proclaim{Theorem 1} There exists a process $(f\circ T^i)$ such that 
$$
  \|E(S_n(f)|\Cal F_0)\|_2 = O\Big(\frac{\sqrt n}{\log^{1/2}n}\Big), \quad \sigma_n^2/n\to 1, \tag4  
$$
but for two different subsequences $(n_k')$, $(n_k'')$, the distributions of $S_{n_k'}/\sqrt{n_k'}$ and
$S_{n_k''}/\sqrt{n_k''}$ converge to different limits.
\endproclaim

\proclaim{Theorem 2} There exists a process $(f\circ T^i)$ such that 
$$
  \|E(S_n(f)|\Cal F_0)\|_2 = O\Big(\frac{\sqrt n}{\log n}\Big), 
$$
but for two different subsequences $(n_k')$, $(n_k'')$, the distributions of $S_{n_k'}/\sigma_{n_k'}$ and
$S_{n_k''}/\sigma_{n_k''}$ converge to different limits.
\endproclaim

\proclaim{Theorem 3} For any sequence of positive reals $c_n\to 0$ there exists a process $(f\circ T^i)$ such that 
$$
  \sum_{n=1}^\infty c_n\frac{\|E(S_n(f)|\Cal F_0)\|_2}{n^{3/2}} < \infty \tag3
$$
but for two different subsequences $(n_k')$, $(n_k'')$, the distributions of $S_{n_k'}/\sigma_{n_k'}$ and
$S_{n_k''}/\sigma_{n_k''}$ converge to different limits.
\endproclaim

From the construction it follows that in Theorems 2,3, the variances $\sigma_n^2$ grow faster than linearly.
It thus remains an open problem whether with a supplementary assumption $\sigma_n^2/n\to const.$, (2) with
$1/2<q<1$ implies the CLT. Also, it remains an open question whether the CLT would hold for $f\in L^{2+\delta}$
for some $\delta>0$.
\bigskip

\subheading{2. Proofs}

\demo{Proof of Theorem 1}
Let $k_n, a_n$, $K_{n,j}$, $n=0,1,\dots$, $j=1,\dots,k_n$, be positive integers, $$
 k_{n+1}=a_nk_n;
$$
$e_{n,j}$, $n=0,1,\dots$, $j=1,\dots,k_n$,
are random variables with $E(e_{n,j})=0$ and $E(e_{n,j}^2)\leq 1$; $k_0=1$, $K_{0,1}=1$,
and $e_{0,1}=0$. We
suppose that $U^ie_{n,j}$, $n=0,1,\dots$, $j=1,\dots,k_n$, $i\in\Bbb Z$ are independent
and 
$$
 K_{n,j} = 3^{j+\sum_{i=1}^{n-1}k_i},\quad \|e_{n,j}\|_2^2 = \frac1{k_n},\quad
j=1,\dots,k_n.
$$
\comment
for $j\leq k_n-1$, $K_{n,j+1}
\geq 2K_{n,j}$ and $K_{n+1,1}\geq 2K_{n,k_n}$. We have 
$$
 \|e_{n,1}\|_2=\dots=\|e_{n,k_n}\|_2
$$ and \endcomment
For $n\geq 1$ even the random variables $e_{n,j}$, $j=1,\dots,k_n$, have values $\pm \sqrt{K_{n,k_n}}$
each with probability $1/(2k_nK_{n,k_n})$ and are zero on the rest of $\Omega$. For $n$ odd $e_{n,j}$, 
$j=1,\dots,k_n$, have the normal distribution $N(0,1/(k_nK_{n,k_n})$.
Notice that
$$
 \|e_{n,j}\|_2^2 = a_n \|e_{n+1,j}\|_2^2.
$$ Define
$$  f = \sum_{n=0}^\infty \sum_{j=1}^{k_n} \sum_{i=0}^{K_{n,j}-1} \frac1{K_{n,j}} U^{-i}
 \Big(-e_{n,j}+\sum_{l=1}^{a_n} e_{n+1,a_n(j-1)+l}\Big).
$$
Then
$$
 \|f\|_2^2 = \sum_{n=0}^\infty \sum_{j=1}^{k_n} \sum_{i=0}^{K_{n,j}-1}
\frac{\|e_{n,j}\|_2^2 +\sum_{l=1}^{a_n}
 \|e_{n+1,a_n(j-1)+l}\|_2^2}{K_{n,j}^2} = 2 \sum_{n=0}^\infty \sum_{j=1}^{k_n}
\frac{\|e_{n,j}\|_2^2}{K_{n,j}}.
$$
which because of exponential growth of the $K_{n,j}$ is finite.

For a positive integer $N$, let $S_N(f) = \sum_{k=0}^{N-1} U^kf$. Because $f$ is regular,
$$
  S_N(f) = \sum_{k=0}^{N-1} U^k \sum_{i\leq 0} P_if = \sum_{k=0}^{N-1} \sum_{i\leq 0}
  P_{i+k} U^kf = S'_N(f) + S''_N(f)
$$
where
$$
 S'_N(f) = \sum_{i=1}^{N-1}\sum_{k=i}^{N-1} P_iU^kf,\quad S''_N(f) = \sum_{i\leq
0}\sum_{k=0}^{N-1} P_iU^kf.
$$
Notice that $S''_N(f) = E(S_N(f)|\Cal F_0)$ and $S'_N(f) = S_N(f) - E(S_N(f)|\Cal F_0)$.

First, we approximate $S'_N(f)$ by a martingale with stationary increments. Notice that
$$
 S'_N(f) = \sum_{i=1}^{N-1}\sum_{k=i}^{N-1} U^iP_0U^{k-i}f =
 \sum_{i=1}^{N-1}\sum_{k=0}^{N-i-1} U^iP_0U^kf.
$$
Let us study the sums $\sum_{k=0}^L P_0U^kf$.\newline
For an integer $K_{n,j}\leq L$ we have
$$\multline
 \sum_{k=0}^L P_0U^k \sum_{i=0}^{K_{n,j}-1} U^{-i}\Big(-e_{n,j}+\sum_{l=1}^{a_n}
 e_{n+1,a_n(j-1)+l}\Big) = \\ 
 = K_{n,j}\Big(-e_{n,j}+\sum_{l=1}^{a_n} e_{n+1,a_n(j-1)+l}\Big).\endmultline \tag1
$$
\comment
because for $k>K_{n,j}-1$ 
$$
 P_0U^k\sum_{i=0}^{K_{n,j}-1} U^{-i}\Big(-e_{n,j}+\sum_{l=1}^{a_n}
e_{n+1,a_n(j-1)+l}\Big) = 0
$$
and for $0\leq k\leq K_{n,j}-1$,
$$
 P_0U^k\sum_{i=0}^{K_{n,j}-1} U^{-i}\Big(-e_{n,j}+\sum_{l=1}^{a_n}
e_{n+1,a_n(j-1)+l}\Big) =
 -e_{n,j}+\sum_{l=1}^{a_n} e_{n+1,a_n(j-1)+l}.
$$
\endcomment
For $1\leq L \leq K_{n,j}$ we have
$$
 \sum_{k=0}^L P_0U^k\sum_{i=0}^{K_{n,j}-1} U^{-i}\Big(-e_{n,j}+\sum_{l=1}^{a_n}
e_{n+1,a_n(j-1)+l}\Big) =
 L \Big(-e_{n,j}+\sum_{l=1}^{a_n} e_{n+1,a_n(j-1)+l}\Big).
$$
For a positive integer $L$ define $K_{n,j}(L) = K_{n,j} = \min\{K_{u,v} : K_{u,v}\geq
L\}$. Then
$$
 \sum_{k=0}^L P_0U^kf = A + B + C
$$
where
$$\multline
 A = \sum_{k=0}^L \sum_{u=1}^{n-1} \sum_{v=1}^{k_n} \sum_{i=0}^{K_{u,v}-1} P_0U^k
\frac1{K_{u,v}} U^{-i}\Big(-e_{u,v}+
 \sum_{l=1}^{a_n} e_{u+1,a_n(v-1)+l}\Big) +\\  + \sum_{k=0}^L \sum_{v=1}^{j-1}
\sum_{i=0}^{K_{n,v}-1} P_0U^k \frac1{K_{n,v}}
 U^{-i}\Big(-e_{n,v} + \sum_{l=1}^{a_n} e_{n+1,a_n(v-1)+l}\Big) = \\
 = \sum_{u=1}^{n-1} \sum_{v=1}^{k_n} \Big(-e_{u,v}+ \sum_{l=1}^{a_n}
e_{u+1,a_n(v-1)+l}\Big) +
 \sum_{v=1}^{j-1} \Big(-e_{n,v} + \sum_{l=1}^{a_n} e_{n+1,a_n(v-1)+l}\Big) =\\
 = \sum_{v=1}^{a_n(j-1)} e_{n+1,v} + \sum_{v=j}^{k_n} e_{n,v},
 \endmultline
$$  $$
 \multline
 B = \sum_{k=0}^L \sum_{i=0}^{K_{n,j}-1} P_0U^k \frac1{K_{n,j}} U^{-i}\Big(-e_{n,j}+
 \sum_{l=1}^{a_n} e_{n+1,a_n(v-1)+l}\Big) = \\
 = \frac{L}{K_{n,j}} \Big(-e_{n,j} + \sum_{l=1}^{a_n} e_{n+1,a_n(v-1)+l}\Big),
 \endmultline
$$
$$
 \multline
 C = \sum_{k=0}^L \sum_{(u,v)\succ (n,j)}  \sum_{i=0}^{K_{u,v}-1} P_0U^k \frac1{K_{u,v}}
U^{-i}\Big(-e_{u,v}+
 \sum_{l=1}^{a_n} e_{u+1,a_n(v-1)+l}\Big) = \\
 = \sum_{(u,v)\succ (n,j)} \frac{L}{K_{u,v}} \Big(-e_{u,v}+ \sum_{l=1}^{a_n}
e_{u+1,a_n(v-1)+l}\Big)
 \endmultline
$$
where $(u,v)\succ (n,j)$ means that either $u>n$ or $u=n$ and $v>j$.\newline
Therefore
$$
 \sum_{k=0}^L P_0U^kf = f'(L) + f''(L)\quad\text{where}\quad f'(L)=A,\,\,\,\,f''(L)=B+C,
$$
hence
$$
 S'_N(f) = \sum_{i=1}^{N-1}U^i f'(N-i-1) + \sum_{i=1}^{N-1}U^if''(N-i-1) = s_1(N)+s_2(N).
$$
Notice that for all $1\leq L\leq N-1$, $\|f'(L)\|_2=1$ and $\|s_1(N)\|_2^2= N-1$.

Suppose that $N = K_{n+1,1}$. Then for $1\leq i\leq K_{n+1,1} - K_{n,k_n}-1$ we have that
$$
 f'(N-i-1) = \sum_{j=1}^{k_{n+1}} e_{n+1,j}.
$$
By the definition, $(1/N)(K_{n+1,1} - K_{n,k_n}) = 1 - 1/3 = 2/3$.

Because $\|e_{n,j}\|_2^2 = 1/k_n$, $\|B\|_2^2 \leq 2/k_n$ and because the $K_{u,v}$
are exponentially increasing, $\|C\|_2^2  = (1/2)\|B\|_2^2$. We deduce that $\|f''(L)\|_2 \to 0$ as
$L\to \infty$ hence
$$
 \|s_2(N)\|_2^2/N \to 0.
$$
Therefore, $S'_N(f)$ can be approximated by the martingale $(s_1(N))_N$ with increments of norm 1,
which for $1\leq i\leq K_{n+1,1} - K_{n,k_n}-1$ have stationary increments of the form 
$U^l\sum_{j=1}^{k_{n+1}} e_{n+1,j}$.
\bigskip

Let us estimate $\|S_N''(f)\|_2$. Note that 
$$
 f = \sum_{n=0}^\infty \sum_{j=1}^{k_n} \sum_{l=1}^{a_n} \sum_{i=0}^{K_{n,j}-1} 
\frac1{K_{n,j}} U^{-i} e_{n+1,a_n(j-1)+l} -
 \sum_{n=0}^\infty \sum_{j=1}^{k_n} \sum_{i=0}^{K_{n,j}-1} \frac1{K_{n,j}} U^{-i} e_{n,j}.
$$

First, let us do several auxilliary calculations.\newline
For $N\geq K$ and a random variable $e\in L^2(\Cal F_0)\ominus L^2(\Cal F_{-1})$ we have
$$
 E\Big(\sum_{l=0}^{N-1} U^l \frac1{K} \sum_{i=0}^{K-1} U^{-i}e|\Cal F_0\Big) =
 \frac1{K} \sum_{i=0}^{K-1} \sum_{l=0}^i U^{l-i}e =  \sum_{i=0}^{K-1} \frac{K-i}{K}U^{-i}e
$$
hence
$$
 \Big\|\sum_{l=0}^{N-1} U^l \frac1{K} \sum_{i=0}^{K-1} U^{-i}e|\Cal F_0\Big\|_2^2 \leq
 K \|e\|_2^2. \tag2
$$
For $N<K$,
$$
 E\Big(\sum_{l=0}^{N-1} U^l \frac1{K} \sum_{i=0}^{K-1} U^{-i}e|\Cal F_0\Big) =
 \frac{N}{K} \sum_{i=0}^{K-N}U^{-i}e + \sum_{i=K-N+1}^{K-1} \frac{K-i}{K} U^{-i}e
$$
hence
$$
 \Big\|\sum_{l=0}^{N-1} U^l \frac1{K} \sum_{i=0}^{K-1} U^{-i}e|\Cal F_0\Big\|_2^2 \leq
 \frac{N^2}{K}\|e\|_2^2. \tag3
$$
In order to estimate $\|S_N''(f)\|_2^2 = \|E(S_N(f)|\Cal F_0)\|_2^2$ we estimate all \newline
$E(\|(\sum_{i=0}^{K_{n,j}-1} \frac1{K_{n,j}} U^{-i} e_{n+1,(a_n-1)j+l}|\Cal F_0)\|_2^2)$
and $E(\|(\sum_{i=0}^{K_{n,j}-1} \frac1{K_{n,j}} U^{-i} e_{n,j}|\Cal F_0)\|_2^2)$.

Let $K_{n,j}(N) = K_{n,j} = \min\{K_{u,v} : K_{u,v}\geq N\}$. The sum $E(S_N(f)|\Cal F_0)$
can be decomposed into summands  $$
 \sum_{i=0}^{K_{u,v}-1} \frac{K_{u,v}-i}{K_{u,v}} U^{-i} e_{u+1,(a_u-1)v+l},\quad
 \sum_{i=0}^{K_{u,v}-1} \frac{K_{u,v}-i}{K_{u,v}} U^{-i} e_{u,v}
$$
for $(n,j)\succ (u,v)$, and $$\gather
 \frac{N}{K_{u,v}} \sum_{i=0}^{K_{u,v}-N}U^{-i}e_{u+1,(a_u-1)v+l} + 
\sum_{i=K_{u,v}-N+1}^{K_{u,v}-1} \frac{K_{u,v}-i}{K_{u,v}} U^{-i}e_{u+1,(a_u-1)v+l}, \\
 \frac{N}{K_{u,v}} \sum_{i=0}^{K_{u,v}-N}U^{-i}e_{u,v} +
 \sum_{i=K_{u,v}-N+1}^{K_{u,v}-1} \frac{K_{u,v}-i}{K_{u,v}} U^{-i}e_{u,v}
 \endgather
$$
if $(u,v)\succ (n,j)$ or $(u,v) = (n,j)$.

By (2) and (3) and exponential growth of $K_{u,v}$ we get that 
$$
  \|S_N''(f)\|_2^2 \approx N\|e_{n,j}\|_2^2 = N/k_n.
$$
We have $K_{n,j-1} < N \leq K_{n,j}$ if $j\geq 2$ and $K_{n-1,k_{n-1}}<N\leq K_{n,j}$ if $j=1$.
Therefore,
$$
  3^{j-1+\sum_{i=1}^{n-1}k_i}<N \leq 3^{j+\sum_{i=1}^{n-1}k_i}. \tag{4}
$$
For $k_n=2^n$ we thus have $k_n = O(\log N)$, therefore
$$
  \|E(S_N(f)|\Cal F_0)\|_2^2 = O(N/\log N).
$$
  
In (4), the number $k_n$ need not depend on $N$.
By an adequate choice of the numbers $k_n$ we can have an increasing sequence $(N_j)$ of integers $N$ such that 
$\|E(S_{N_j}(f)|\Cal F_0)\|_2 \to 0$ as $j\to\infty$.
\bigskip

We proved that $\|E(S_N(f)|\Cal F_0)\|_2/\sqrt N\to 0$, $\|s_2(N)\|_2/\sqrt N\to 0$, and $\|s_1(N)\|_2 = \sqrt{N-1}$.
Therefore, the limit behaviour of $S_N(f)/\sqrt N$ will be the same as that of $s_1(N)/\sqrt N$.
The term $s_1(N)$ is a sum of martingale differences (independent random variables, in fact) $U^if'(N-i-1)$, 
each of the norm $\|U^if'(N-i-1)\|_2=1$,
$1\leq i\leq N-1$. Let us consider a subsequence $(N_n)$, $N_n = K_{n,1}$. Then for $n$ even, $2/3$ of the summands in
$s_1(N)$ will have a symmetrized Poisson distribution while for $n$ odd they will be normally distributed.
\enddemo
\qed
\bigskip

\demo{Proof of Theorems 2 and 3} 
Let $n_k = 2^k$,
$e_k$ random variables with
$\|e_k\|_2 = a_k/k$, $0\leq a_k\leq 1$, $\sum_{k=1}^\infty a_k/k =\infty$, such that for each $k$, 
$U^ie_k$ are independent, and let
$$
  f = \sum_{k=1}^\infty \frac1{n_k} \sum_{i=0}^{n_k-1} U^{-i}e_k.
$$
We have $\|f\|_2 \leq \sum_{k=1}^\infty \|e_k\|_2/\sqrt{n_k} <\infty$ due to the exponential growth of the $n_k$s.

For a positive integer $N$ we have
$$
  S_N(f) = \sum_{k=1}^\infty \sum_{j=0}^{N-1} \sum_{i=0}^{n_k-1} \frac1{n_k} U^{j-i}e_k = S'_N(f) + S''_N(f)
$$
where
$$
  S'_N(f) = S_N(f) - E(S_N(f)|\Cal F_0) = \sum_{k=1}^\infty \sum_{j=0}^{N-1} \sum_{i=0}^{(j\wedge n_k)-1} 
  \frac1{n_k} U^{j-i}e_k
$$
($j\wedge n_k = \min\{j, n_k\}$) and
$$
  S''_N(f) = E(S_N(f)|\Cal F_0) = \sum_{k=1}^\infty \sum_{j=0}^{N-1} \sum_{i=j}^{n_k-1} \frac1{n_k} U^{j-i}e_k.
$$

First, let us study $S'_N(f)$. \newline
For $N\leq n_k$ we have
$$
  \sum_{j=0}^{N-1} \sum_{i=0}^{(j\wedge n_k)-1} U^{j-i}e_k = \sum_{j=1}^{N-1} (N-j)U^je_k \tag{5}
$$
and for $N>n_k$ we have
$$
  \sum_{j=0}^{N-1} \sum_{i=0}^{(j\wedge n_k)-1} U^{j-i}e_k = 
  \sum_{j=1}^{N-n_k} n_kU^je_k + \sum_{j=N-n_k+1}^{N-1} (N-j)U^je_k. \tag{6}
$$
For $l\geq N$ and $l\leq 0$ we have $P_lS_N(f) = 0$. For $1\leq l\leq N-1$ we thus, using (5) and (6), get
$$
  P_lS_N(f) = \sum_{k\geq 1: n_k\leq N-l} U^je_k + \sum_{k\geq 1: n_k\geq N+1-l} \frac{N-l}{n_k} U^le_k. \tag{7}
$$
For any $1>\epsilon>0$ (arbitrarilly small) and $1\leq l\leq N-\epsilon N=N(1-\epsilon)$,  
$$
  \sum_{k\geq 1: n_k\leq N-l} U^le_k = \sum_{k\geq 1: n_k\leq \epsilon N} U^le_k + 
  \sum_{k\geq 1: \epsilon N<n_k\leq N} U^le_k.
$$
Let $[x]$ denoter the integer part of $x$. We have
$$
  S'_N(f) = \sum_{l=1}^{N-1} P_lS_N(f) = \sum_{l=1}^{[N(1-\epsilon)]} P_lS_N(f) +
  \sum_{l=[N(1-\epsilon)]+1}^{N} P_lS_N(f)
$$
where
$$\multline
  \sum_{l=1}^{[N(1-\epsilon)]} P_lS_N(f) = \\
  \sum_{l=1}^{[N(1-\epsilon)]} U^l \sum_{k\geq 1: n_k\leq N-l} e_k +
  \sum_{l=1}^{[N(1-\epsilon)]} U^l \sum_{k\geq 1: n_k\geq N+1-l} \frac{N-l}{n_k} e_k = \\
  = \sum_{l=1}^{[N(1-\epsilon)]} U^l \sum_{k\geq 1: n_k\leq \epsilon N} e_k +
  \sum_{l=1}^{[N(1-\epsilon)]} U^l \sum_{k\geq 1: \epsilon N<n_k\leq N} e_k + \\
  + \sum_{l=1}^{[N(1-\epsilon)]} U^l \sum_{k\geq 1: n_k\geq N+1-l} \frac{N-l}{n_k} e_k.
  \endmultline
$$
Because $n_k=2^k$, 
$$
  \Big\|\sum_{k\geq 1: n_k\geq N+1-l} \frac{N-l}{n_k} U^le_k\|_2 \leq 2\|e_k\Big\|_2 \leq 2/k;
$$ 
$\epsilon N <n_k\leq N$ if and only if $\log N + \log(1-\epsilon)<k\leq \log N$. We thus deduce, for $b(N) = 
\|\sum_{k=1}^{\log N} e_k\|_2$ ,
$$
  \lim_{N\to\infty} \frac1{b(N)\sqrt N} \Big\|\sum_{l=1}^{[N(1-\epsilon)]} P_lS_N(f) - \sum_{l=1}^{[N(1-\epsilon)]} U^l
  \sum_{k\geq 1: n_k\leq \epsilon N} e_k\Big\|_2 = 0.
$$
For all $[N(1-\epsilon)]+1\leq l\leq N$ we have, by (7),
$\|P_lS_N(f)\|_2 \leq b(N) + 2/k$ hence 
$$
  \lim_{\epsilon\searrow 0}\lim_{N\to\infty} \frac1{b(N)\sqrt N} \Big\|\sum_{l=[N(1-\epsilon)]+1}^{N-1} P_lS_N(f) - 
  \sum_{l=[N(1-\epsilon)]+1}^{N-1} \sum_{k\geq 1: n_k\leq N(1-\epsilon)} U^le_k\Big\|_2 = 0
$$
Therefore,
$$\multline
  \lim_{N\to\infty} \frac1{b(N)\sqrt N} \Big\|S'_N(f) - \sum_{l=0}^{N-1} U^l \sum_{k\geq 1: n_k\leq N} e_k\Big\|_2 = \\
  \lim_{\epsilon\searrow 0}\lim_{N\to\infty} \frac1{b(N)\sqrt N} \Big\|\sum_{l=1}^{N-1} P_lS_N(f) - 
  \sum_{l=1}^{N-1} \sum_{k\geq 1: n_k\leq N(1-\epsilon)} U^le_k\Big\|_2 = 0.
  \endmultline
$$

Now, let us study $S''_N(f) = E(S_N(f)|\Cal F_0)$.\newline
For $N\leq n_k$ we have
$$
  \sum_{j=0}^{N-1} \sum_{i=j}^{n_k-1} U^{j-i}e_k = \sum_{i=0}^{n_k-N} N U^{-i}e_k + \sum_{i=n_k-N+1}^{n_k-1}
  (n_k-i) U^{-i}e_k
$$
and for $N>n_k$ we have
$$
  \sum_{j=0}^{N-1} \sum_{i=j}^{n_k-1} U^{j-i}e_k = \sum_{j=0}^{n_k-1} (n_k-j) U^{-j}e_k,
$$
hence
$$\multline
  S''_N(f) = \sum_{k\geq 1: n_k< N} \sum_{j=0}^{n_k-1} \frac{n_k-j}{n_k} U^{-j}e_k + \\
  + \sum_{k\geq 1: n_k\geq N} \Big[\sum_{j=0}^{n_k-N} \frac{N}{n_k} U^{-j}e_k +
  \sum_{j=n_k-N+1}^{n_k-1} \frac{n_k-j}{n_k} U^{-j}e_k \Big].
  \endmultline
$$
We have
$$\gather
  (1/\sqrt 6)\|e_k\|_2\sqrt{n_k}\leq \Big\|\sum_{j=0}^{n_k-1} \frac{n_k-j}{n_k} U^{-j}e_k\Big\|_2 \leq \|e_k\|_2\sqrt{n_k},
  \quad n_k<N,\\
  \Big\|\sum_{j=0}^{n_k-N} \frac{N}{n_k} U^{-j}e_k + \sum_{j=n_k-N+1}^{n_k-1} \frac{n_k-j}{n_k} U^{-j}e_k \Big\|_2
  \leq \frac{N}{\sqrt{n_k}} \|e_k\|_2, \quad n_k\geq N.
  \endgather
$$
  
Because $n_k = 2^k$ grow exponentially fast, for a constant $0<c<\infty$ not depending on $N$ and $[x]$ the integer part 
of $x$ we get $\sum_{k\geq 1: n_k\geq N} (N/\sqrt{n_k})\|e_k\|_2 \leq c\|e_{[\log N]}\|_2$.
We deduce that for some constants $c',c''>0$ we have
$$
  c'N \|e_{[\log N]}\|_2^2< \|E(S_N(f)|\Cal F_0)\|_2^2 <c''N \|e_{[\log N]}\|_2^2.
$$
Because $\|e_k\|_2 = a_k/k$,
$$
  c'\frac{Na^2_{[\log N]}}{\log^2 N} < \|E(S_N(f)|\Cal F_0)\|_2^2 <c''\frac{Na^2_{[\log N]}}{\log^2 N}.\tag8
$$
\medskip
From (8) we deduce that the series $\sum_{n=1}^\infty n^{-3/2}\|E(S_n(f)|\Cal F_0)\|_2$ converges if and only if
$\sum_{n=1}^\infty a_{[\log n]}/(n[\log n])$ converges; this is equivalent to the convergence of
$\sum_{n=1}^\infty a_n/n$.
If $c_n$ are positive real numbers, $c_n\to 0$, we can choose the $a_n$ so that
$$
  \sum_{n=1}^\infty c_n\frac{\|E(S_n(f)|\Cal F_0)\|_2}{n^{3/2}} \approx 
  \sum_{n=1}^\infty \frac{c_na_{[\log n]}}{n[\log n]} <\infty
$$
but $\sum_{n=1}^\infty a_n/n = \infty$.
Then we can choose $e_n$ so that $S_N(f)/\|S_N(f)\|_2$ converge along subsequences to different limit points.
\enddemo
\qed


\Refs
\widestnumber\key{Wu-Wo}
\ref \key B 1961\by Billingsley, P. \paper The Lindeberg-L\'evy theorem for martingales \jour Proc. Amer. Math. Soc.
\vol 12 \pages 788-792 \yr 1961 \endref 
\ref \key D-L  \by Derriennic, Y and Lin, M. \paper The central limit theorem for Markov
chains with normal transition operators, started at a point \jour Probab. Theory Relat. Fields
\vol 119 \pages 509-528 \yr 2001 \endref
\ref \key G \by Gordin, M.I. \paper A central limit theorem for stationary processes \jour Soviet Math. Dokl. \vol 10 
\pages 1174-1176 \yr 1969 \endref
\ref \key G-Ho \by Gordin, M.I. and Holzmann, H. \paper The central limit theorem for stationary Markov chains under invariant splittings 
\jour Stochastics and Dynamics \vol 4 \pages 15-30 \yr 2004 \endref 
\ref \key G-L \by Gordin, M.I. and Lif\v sic, B.A. \paper A remark about a Markov process with normal transition operator
\paperinfo In: Third Vilnius Conference on Probability and Statistics \vol 1 \pages 147-148 \yr 1981 \endref
\ref \key Ha-He \by Hall, P. and Heyde, C.C. \book Martingale Limit Theory           
and its Application \publ Academic Press \publaddr New York \yr 1980 \endref 
\ref \key I \by Ibragimov, I.A.  \paper A central limit theorem for a class of dependent random variables 
\jour Theory Probab. Appl. \vol 8 \pages 83-89 \yr 1963 \endref
\ref \key K-V \by Kipnis, C. and Varadhan, S.R.S. \paper Central limit theorem for additive functionals of reversible
Markov processes and applications to simple exclusions \jour Comm. Math. Phys. \vol 104 \pages 1-19 \yr 1986 \endref
\ref \key Kl-Vo 1 \by Klicnarov\'a, J. and Voln\'y, D. \paper An invariance principle for non adapted processes
\paperinfo preprint \yr 2007 \endref
\ref \key Kl-Vo 2 \by Klicnarov\'a, J. and Voln\'y, D. \paper Exactness of a Wu-Woodroofe's approximation with linear growth of variances
\paperinfo preprint \yr 2007 \endref
\ref \key M-Wo \by Maxwell, M. and Woodroofe, M. \paper Central limit theorems for additive
functionals of Markov chains \jour Ann. Probab. \vol 28 \pages 713-724 \yr 2000 \endref
\ref \key P-U \by Peligrad, M. and Utev, S. \paper A new maximal inequality
and invariance principle for stationary sequences \jour Ann. Probab.
\vol 33 \pages 798-815 \yr 2005 \endref
\ref \key V 2006 \by Voln\'y, D. \paper Martingale approximation of non adapted stochastic processes
with nonlinear growth of variance \paperinfo Dependence in Probability and Statistics
Series: Lecture Notes in Statistics, Vol. 187 Bertail, Patrice; Doukhan, Paul; Soulier, Philippe (Eds.) 
\yr 2006 \endref
\ref \key Wo \by Woodroofe, M. \paper A central limit theorem for functions of a Markov chain with applications to
shifts \jour Stoch. Proc. and their Appl. \vol 41 \pages 31-42 \yr 1992 \endref
\ref  \key Wu-Wo \by Wu, Wei Biao and Woodroofe, M. \paper Martingale approximation for 
sums of stationary processes \jour Ann. Probab. \vol 32 \pages 1674-1690 \yr 2004 \endref  
\endRefs

\Refs
\widestnumber\key{Wu-Wo}
\ref \key B 1961\by Billingsley, P. \paper The Lindeberg-L\'evy theorem for martingales \jour Proc. Amer. Math. Soc.
\vol 12 \pages 788-792 \yr 1961 \endref 
\ref \key B 1968\by Billingsley, P. \book Convergence of Probability Measures \publ
John Wiley \publaddr New York \yr 1968 \endref
\ref \key Du-Gol \by D\"urr, D. and Goldstein, S. \paper Remarks on the CLT for weakly dependent random variables
\paperinfo in: S. Albeverio, Ph. Blanchard, and L. Streit, eds., Stochastic Processes - Mathematics and Physics,
Proc. Bielefeld, 1984, Lecture Notes in Maths. No 1158 \endref
\ref \key G \by Gordin, M.I. \paper A central limit theorem for stationary processes \jour Soviet Math. Dokl. \vol 10 
\pages 1174-1176 \yr 1969 \endref
\ref \key Ha-He \by Hall, P. and Heyde, C.C. \book Martingale Limit Theory           
and its Application \publ Academic Press \publaddr New York \yr 1980 \endref 
\ref \key I \by Ibragimov, I.A.  \paper A central limit theorem for a class of dependent random variables 
\jour Theory Probab. Appl. \vol 8 \pages 83-89 \yr 1963 \endref
\ref \key M-Wo \by Maxwell, M. and Woodroofe, M. \paper Central limit theorems for additive
functionals of Markov chains \jour Ann. Probab. \vol 28 \pages 713-724 \yr 2000 \endref
\ref \key V 1987 \by Voln\'y, D. \paper Martingale decompositions of stationary processes 
\jour Yokohama Math. J. \vol 35 \pages 113-121 \yr 1987 \endref
\ref \key V 1993 \by Voln\'y, D. \paper Approximating martingales and the CLT         
for strictly stationary processes \jour Stoch. Proc. and their
Appl. \vol 44 \pages 41-74 \yr 1993 \endref   
\ref \key V 2006 \by Voln\'y, D. \paper Martingale approximation of non adapted stochastic processes
with nonlinear growth of variance \paperinfo Dependence in Probability and Statistics
Series: Lecture Notes in Statistics, Vol. 187 Bertail, Patrice; Doukhan, Paul; Soulier, Philippe (Eds.) 
\yr 2006 \endref
\ref \key Wo \by Woodroofe, M. \paper A central limit theorem for functions of a Markov chain with applications to
shifts \jour Stoch. Proc. and their Appl. \vol 41 \pages 31-42 \yr 1992 \endref
\ref  \key Wu-Wo \by Wu, Wei Biao and Woodroofe, M. \paper Martingale approximation for 
sums of stationary processes \jour Ann. Probab. \vol 32 \pages 1674-1690 \yr 2004 \endref  
\endRefs

\Refs
\ref \by Maxwell and Woodroofe \endref
\ref \by Wu and Woodroofe \endref
\ref \by Peligrad and Utev\endref
\ref \by Klicnarov\'a\endref
\ref \by Voln\'y\endref
\ref \by Hall and Heyde\endref
\endRefs
\enddocument